\documentclass[10pt,reqno]{amsart}
\usepackage{amsmath, amsthm, amssymb, stmaryrd}
\usepackage{hyperref}
\usepackage{enumerate}
\usepackage{url}
\usepackage{color}
\usepackage{tikz} 

\let\OLDthebibliography\thebibliography
\renewcommand\thebibliography[1]{
  \OLDthebibliography{#1}
  \setlength{\parskip}{0pt}
  \setlength{\itemsep}{0pt plus 0.3ex}
}

\usepackage{mathtools}

\usepackage{multirow, array}
\usepackage{placeins}

\usepackage{caption} 
\advance \topmargin by -\headheight
\advance \topmargin by -\headsep
     
\evensidemargin \oddsidemargin
\newtheorem{thm}{Theorem}[section]

\newtheorem{prop}[thm]{Proposition}

\newtheorem{conj}[thm]{Conjecture}

\theoremstyle{definition}
\newtheorem{defn}[thm]{Definition}

\theoremstyle{remark}

\numberwithin{equation}{section}

\newcommand*\wrapletters[1]{\wr@pletters#1\@nil}
\def\wr@pletters#1#2\@nil{#1\allowbreak\if&#2&\else\wr@pletters#2\@nil\fi}

\def\eps{\varepsilon}

\def\le{\leqslant} \def\ge{\geqslant}

\def \bN {\mathbb N}
\def \bP {\mathbf P}

\def \bR {\mathbb R}
\def \bZ {\mathbb Z}

\def \bL {\mathbf L}

\def \bq {\mathbf q}

\def \bx {\mathbf x}

\def \by {\mathbf y}

\def \bbR {\mathbb R}

\def \cC {\mathcal C}
\def \cD {\mathcal D}

\def \dim {\mathrm{dim}}

\def \Bad {{\mathbf{Bad}}}

\newcommand{\DI}{\mathbf{DI}}
\newcommand{\Sing}{\mathbf{Sing}}

\begin{document}
\title{Diophantine sets and Dirichlet improvability}
\author[Antoine Marnat]{Antoine Marnat\\Moscow Center of Fundamental and Applied Mathematics}
\address{Moscow Center of Fundamental and Applied Mathematics}
\email{antoine.marnat@gmail.com}

\subjclass[2010]{}
\thanks{Supported by FWF-Projekt I 5554 "Diophantine approximations, arithmetic sequences \& analytic number theory"}
\date{}
\begin{abstract}
This note pushes further the discussion about relations between Dirichlet improvable, badly approximable and singular points held in \cite{yorkies} by considering Diophantine sets extending the notion of badly approximability.\end{abstract}

\maketitle

\section{Introduction}
The aim of this paper is to extend slightly the main result of \cite{yorkies}, questioning relations between Dirichlet improvability, badly approximability and singularity. We first provide a short description of this setting and definitions.\\

For $x \in \bbR$, we denote $\langle x \rangle := \min\{|x-z| \mid z \in \bZ\}$ the distance from $x$ to a nearest integer. For $\bx,\by \in \bR^n$, we denote $\bx\cdot \by := x_1y_1 + \cdots + x_ny_n$ the usual scalar product.\\

The celebrated \emph{Dirichlet's Theorem}, root of Diophantine approximation, reads as follows.

\begin{thm}[Dirichlet, 1842]
Let $\bx\in\bR^n$. For every $Q>1$ there exists $\bq\in \bZ^n$ such that
\begin{equation*}
0< |\bq| \le Q  \qquad \textrm{ and } \qquad \langle \bx\cdot \bq \rangle \le Q^{-n}.
\end{equation*}
\end{thm}

A main interest in Diophantine approximation is to question when and how Dirichlet's theorem can be improved. This lead Davenport and Schmidt \cite{DS} to call a given $\bx\in\bR^n$ \emph{Dirichlet improvable} if there exists $\eps \in (0,1)$ such that for all sufficiently large $Q$, there exists $\bq$ with
\begin{equation}\label{DI}\tag{DI}
0< |\bq| \le Q  \qquad \textrm{ and } \qquad \langle \bx\cdot \bq \rangle \le \eps Q^{-n}.\end{equation}
We denote by $\DI_n(\eps)$ the set of $\bx$ satisfying \eqref{DI}, so that the set of Dirichlet improvable numbers is 
\[ \DI_n = \bigcup_{\eps \in (0,1)} \DI_n(\eps).\]
Furthermore, $\bx \in \bR^n$ is called \emph{singular} if it is in $\DI_n(\eps)$ for $\eps$ arbitrarily small. Denote \[\Sing_n := \bigcap_{\eps \in (0,1)} \DI_n(\eps).\]
On the opposite, we say that $\bx\in\bR^n$ is \emph{badly approximable} if there exists $\eps = \eps(\bx) \in (0,1)$ such that
\[ \langle \bq\cdot \bx \rangle \ge \eps |\bq|^{-n}\]
for all $\bq\in\bR^n$. We denote the set of badly approximable numbers $\Bad_n$.\\

A very natural question is the size and relation between these sets. Regarding sizes, it is well known that $\DI_n$, $\Bad_n$ and $\Sing_n$ have $0$ Lebesgue measure. The set $\Bad_n$ has full Hausdorff dimension, hence $\DI_n$ as well: 
\[ \dim_H(\Bad_n) = \dim_H(\DI_n) = n.\]
For $\Sing_m$, the study of the Hausdorff dimension is much more involved. The result requires the powerful variational principle in parametric geometry of numbers by Das, Fishman, Simmons and Urba\'nski \cite{David}, extending to a wider setting (including ours: approximation to a linear form) a brillant result by Cheung and Chevallier \cite{CC} for simultaneous approximation.
\[\dim_H(\Sing_n) = \frac{n^2}{n+1}, \qquad \textrm{ for } n\ge 2.\]

Regarding relations, we have the inclusion of the disjoint union $\Bad_n \sqcup \Sing_n \subset \DI_n$. When $n=1$, singular reals are the rationals and the last inclusion is an equality. In higher dimension, one can ask about the set $$\mathcal{FS}_n:= \DI_n \setminus \left( \Bad_n \sqcup \Sing_n\right).$$ This set is studied in \cite{yorkies}, where it is proved to be uncountable. However, it is probably far from the truth and one conjectures

\begin{conj}
$\dim_H\left( \mathcal{FS}_n\right) = n.$
\end{conj}

Actually, in \cite{yorkies} the result is finer than just uncountability. This requires the definition of exponents of Diophantine approximation. Namely, for $\bx\in\bR^m$, we define $\omega(\bx)$ as the supremum of positive reals $w$ such that
\[0 < |\bq| \le Q \; , \qquad \langle \bq\cdot \bx \rangle \le Q^{-w}  \]
has integer solution $\bq$ for arbitrarily large $Q$.\\
This exponent is usually referred to as exponent of Diophantine approximation to a linear form. It can take any value in the interval $[n,\infty]$. The main result in \cite{yorkies} asserts that there are uncountably many $\bx \in \bR^n$ in 
\[ \left(\DI_n(\eps) \setminus \DI_n(c_n\eps)\right) \setminus \left(\Bad_n \sqcup \Sing_n\right)\]
with prescribed exponent $\omega(\bx) \in [n,\infty]$. Here $c_n$ is an explicit constant. The exclusion of $\Sing_n$ is already included in the exclusion of $\DI_n(c_n\eps)$, but we express this set to enlighten its link to $\mathcal{FS}_n$.\\
The motivation of this paper is the obvious observation that $\bx \notin \Bad_n$ if $\omega(\bx) > n$. This leads to consider the following notion : fix dimension $n$ and let $\eps \in (0,1)$ and finite $w \ge n$. We define the $(\eps,w)$-Diophantine set
\[ \cD_{w}(\eps) := \left\{ \bx\in \bR^n \mid \langle \bq\cdot \bx \rangle \ge \eps |\bq|^{-w} >0, \quad \textrm{ for all } \bq \in \bN^n  \right\}\]
and denote $\cD_{w} := \bigcup_{\eps \in (0,1)} \cD_{w}(\eps)$. One can see that $\cD_{n} = \Bad_n$.\\

These Diophantine sets play a role in dynamical systems. For example they describe the Diophantine condition in
small divisors problems with applications to KAM theory, Aubry-Mather theory,
conjugation of circle diffeomorphisms, and so on (see \cite{KAM1,KAM2,KAM3,KAM4,KAM5}).\\

Our main theorem reads as follows.\\

\begin{thm}\label{ThmMain}
Fix the dimension $n\ge2$, a finite exponent $w > n$, parameters $\eps,\nu \in (0,1)$ and the constants $c_n = e^{-20(n+1)^3(n+10)}$ and $c_n'= e^{-20(n+1)^2(w+1)(n+10)}$. Then
\begin{itemize}
\item 
there exists uncountably many $$\bx \in \left( \DI_n(\eps) \setminus \DI_n(c_n\eps) \right) \setminus \cD_{w}$$
with $\omega(\bx) = w$,
\item there exists uncountably many $$\bx \in \left( \DI_n(\eps) \setminus \DI_n(c_n\eps) \right) \cap \left( \cD_{w}(\nu) \setminus \cD_{w}(c_n'\nu) \right).$$
\end{itemize}
\end{thm} 

As discussed in \cite{yorkies}, the constants $c_n$ and $c_n'$ are not optimized.\\

So far, we considered only approximation to a linear form. An analogous argumentation applies for simultaneous approximation, where we replace \eqref{DI} by
\[ 0 < q < Q \textrm{ and } \|q\bx\| \le \eps Q^{-1/m}\]
and $\cD_w(\eps)$  by
\[ \cD_{\lambda}(\eps) := \left\{  \bx \in \bR^m \mid \|q\bx\| \ge \eps |q|^{-\lambda}>0 \right\}\]
where $\|.\|$ is the distance to a nearest integer point. We get uncountably many points $\bx \in \bR^m$ in 
\[\bx \in \left( \DI_n(\eps) \setminus \DI_n(c_n\eps) \right) \setminus \cD_{\lambda}\]
with exponent of simultaneous approximation $\lambda$ and uncountably many $\bx \in \bR^m$ in the intersection
\[\bx \in \left( \DI_n(\eps) \setminus \DI_n(c_n\eps) \right) \cap \left( \cD_{\lambda}(\nu) \setminus \cD_{\lambda}(c_n'\nu) \right).\]

One could also consider approximation by rational subspaces of dimension exactly $d$, for $0\le d < n$, as introduced by Laurent \cite{ML} following  Schmidt \cite{SchmidtH}. However, their usual study via compounds convex bodies involves constants depending on $\bx$ (See \cite[Proposition 3.1]{RoySpec} or \cite[\S 4]{BL2010}) that seem to break our proof. See \cite{yorkies} for more discussion about the sets $\DI_n$, $\Bad_n$ and $\Sing_n$ related to these approximation settings.\\


The proof of Theorem \ref{ThmMain} relies on the parametric geometry of numbers and Roy's fundamental theorem \cite{Roy}. We provide a short introduction to parametric geometry of numbers in Section \ref{SecPGN}, and prove our main theorem in Section \ref{SecProof}.\\

\section{Parametric geometry of numbers}\label{SecPGN}
Parametric geometry of numbers was developped by Schmidt and Summerer \cite{SS1,SS2}, answering a question of Schmidt \cite{SchLum}. It was pushed by a fundamental theorem of Roy \cite{Roy}, that was quantified and extended to a matrix setting by Das, Fishman, Simmons and Urba\'nski \cite{David}.\\

Fix dimension $n>1$ and $\bx \in \bR^n$. Consider the convex body with parameter $q\ge0$
\[\cC_{\bx}(e^q) := \{ \by\in\bR^{n+1} \mid |y_i| \le 1, |\by\cdot (1,\bx) \le e^{-q} \} \]
and its $d$th successive minima
\[ \lambda_{\bx,d}(q) := \lambda_d (\bZ^{n+1}, \cC_{\bx}(e^q))\]
for any $d$ between $1$ and $n+1$. Following Schmidt and Summerer, we consider the successive minima map consisting of their logs
\[\bL_{\bx} : [0,\infty) \to \bR^{n+1} : q \mapsto \bL_{\bx}(q) := (L_{\bx,1}(q), \ldots , L_{\bx,n+1}(q)),\]
where $L_{\bx,d}(q) := \log \lambda_{\bx,d}(q)$. It appears that the map $\bL_{\bx}$ encodes Diophantine properties of $\bx$.

\begin{prop}\label{GPNDP}
Fix dimension $n$ and $\bx\in\bR^n$.
\begin{itemize}
\item{} $\bx \in \DI_n(\eps)$  if and only if for all sufficiently large $q$
\begin{equation}
\frac{q}{n+1} - L_{\bx,1}(q) \ge -\frac{\log (\eps)}{n+1}
\end{equation}
\item {}$\bx \in \cD_{w}(\eps)$  if and only if for all sufficiently large $q$
\begin{equation}
\frac{q}{w+1} - L_{\bx,1}(q) \le -\frac{\log (\eps)}{w+1}
\end{equation}
\item $\liminf_{q\to \infty} \frac{L_{\bx,1}(q)}{q} = \frac{1}{1+\omega(\bx)}$
\end{itemize}
\end{prop}

See \cite[Proposition 3.1]{RoySpec} and \cite[Lemma 2.1]{yorkies} for proofs.\\

The following notion of \emph{system} was introduced by Roy \cite[Definition 4.5]{RoySpec}. The latter approach exactly the familly of successive minima maps.

\begin{defn}
Be $I$ a subinterval of $[0,\infty]$ with nonempty interior. A system on $I$ is a continuous linear map $\bP=(P_1, \ldots , P_{n+1}) : I \to \bR^{n+1}$ with the following properties.
\begin{enumerate}[(i)]
\item{}For each $q\in I$, $0\le P_1(q) \le P_2(q) \le \cdots \le P_{n+1}(q)$ and $ P_1(q) + P_2(q) + \cdots + P_{n+1}(q) = q$.
\item{}If $I'\subset I$ is a nonempty open subinterval on which $\bP$ is differentiable, then there exists integers $r_1,r_2$ with $1 \le r_1 \le r_2 \le n+1$ such that $P_{r_1}, \ldots , P_{r_2}$ coincide on $I'$ and have slope $\frac{1}{r_2-r_1+1}$ while all other component $P_i$ are constant on $I'$.
\item{} If $q$ is an interior point of $I$ at which $\bP$ is not differentiable, and if $r_1,r_2, s_1,s_2$ are integers such that
$$ P_i'(q^-) = \frac{1}{r_2-r_1+1} \quad (r_1\le i \le r_2) \textrm{ and }  P_j'(q^+) = \frac{1}{s_2-s_1+1} \quad (s_1\le j \le s_2)  $$
and if $r_1 \le s_2$ then $P_{r_1}(q)= \cdots = P_{s_2}(q)$
\end{enumerate}
 
\end{defn}

Roy's fundamental theorem reads as follows \cite[Corollary 4.7]{RoySpec}, \cite[Theorems 1.3 \& 1.8]{Roy}.

\begin{thm}[Roy, 2015]\label{Roy}
Fix dimension $n\ge1$ and $q_0 \ge0$. For each $\bx \in \bR^n$, there exists a system $\bP:[q_0,\infty) \to \bR^{n+1}$ such that $L_{\bx}-\bP$ is bounded on $[q_0,\infty)$. Conversely, for each system $\bP:[q_0,\infty) \to \bR^{n+1}$, there exists a $\bx \in \bR^n$ such that $L_{\bx}-\bP$ is bounded. In particular, for each $q \ge q_0$ and a constant $R_n$
\[| L_{\bx}-\bP | \le R_n.\]
\end{thm}

The constant $R_n$ induces the constants $c_n = e^{4(n+1)R_n}$ and $c_n' = e^{4(w+1)R_n}$ in Theorem \ref{ThmMain}.\\

We call two systems $\bP_1$ and $\bP_2$ non-equivalent if there exists $q$ such that $|\bP_1(q) - \bP_2(q)| > 2C_n$. No point in $\bR^n$ has successive minima map close to two non-equivalent systems in the sense of Theorem \ref{Roy}.

\section{Proof of Theorem \ref{ThmMain}}\label{SecProof}

Fix dimension $n>1$, exponent $w>n$ and $\eps,\nu\in(0,1)$. In view of Theorem \ref{Roy} and Proposition \ref{GPNDP}, the proof of Theorem \ref{ThmMain} reduces to the construction of uncountably many non-equivalent systems satisfying

\begin{itemize}
\item{}$\liminf_{q\to \infty} \frac{P_1(q)}{q} = \frac{1}{1+ w}$,\\
\item{}$\liminf_{q\to \infty} \left(\frac{q}{n+1} - P_1(q)\right) = -\frac{\log (\eps)}{n+1} + 2R_n$,\\
\item{}Either $\limsup_{q\to\infty} \left(\frac{q}{w+1} - P_{1}(q)\right) = (-\frac{\log (\nu)}{w+1} + 2R_n)$ \\
or $\limsup_{q\to\infty} \left(\frac{q}{w+1} - P_{1}(q)\right) = \infty$.
\end{itemize}

Note that the second and third inequalities provide the definition of the constants $c_n := e^{4(n+1)R_n}$ and $c_n' = e^{4(w+1)R_n}$.\\

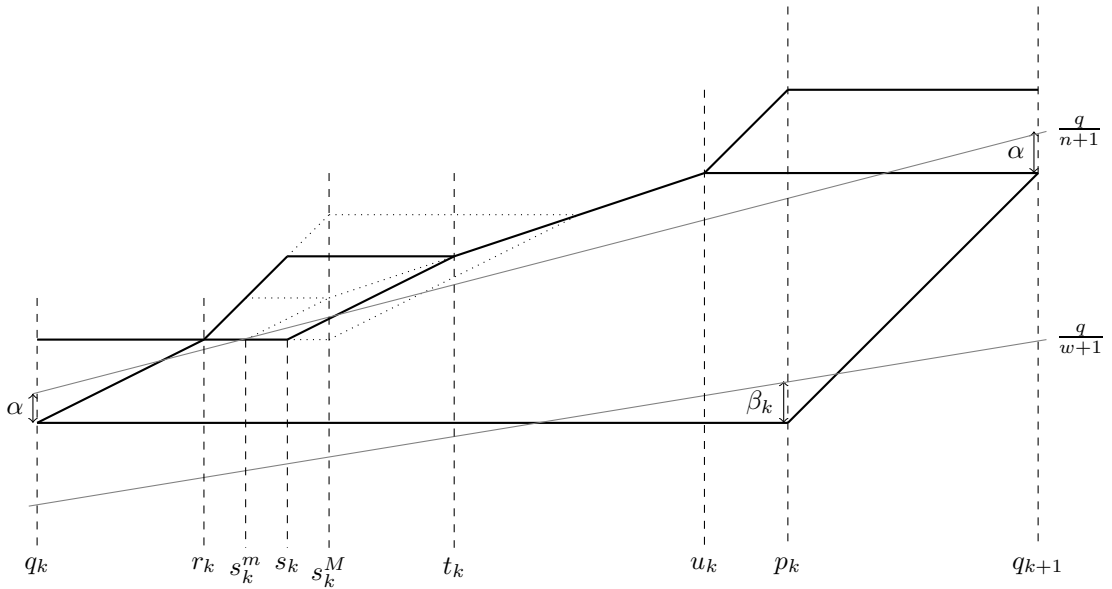
\begin{figure}[h!]
 \begin{center}
 \begin{tikzpicture}[scale=0.55]
 
 \draw[black, thick] (0,-1) -- (4,1) ;
 \draw[black, thick] (0,-1) -- (18,-1) ;
 \draw[black, thick] (0,1) -- (6,1) ;
 \draw[black, thick] (6,3) -- (4,1) ;
 \draw[black, thick] (6,1) -- (10,3) ;
 \draw[black, thick] (6,3) -- (10,3) ;
 \draw[black, thick] (10,3) -- (16,5) ;
 
  \draw[black, dotted] (5,2) -- (7,2) ;  
  \draw[black, dotted] (5,1) -- (7,2) ;
  \draw[black, dotted] (7,2) -- (10,3) ;
  \draw[black, dotted] (6,1) -- (7,1) ;
  \draw[black, dotted] (6,3) -- (7,4) ;
  \draw[black, dotted] (7,4) -- (13,4) ;
  \draw[black, dotted] (7,1) -- (13,4) ;

 \draw[black, thick] (16,5) -- (18,7) ;
 \draw[black, thick] (18,7) -- (24,7) ;
 \draw[black, thick] (24,5) -- (18,-1) ; 
 \draw[black,thick] (24,5)-- (16,5);

 \draw[black, <->] (-0.1,-0.3) -- (-0.1,-1) node [black, midway, left] {$\alpha$};
 \draw[black, <->] (17.9,0) -- (17.9,-1) node [black, midway, left] {$\beta_k$};
  \draw[black, <->] (23.9,5) -- (23.9,6) node [black, midway, left] {$\alpha$};
    
 \draw[gray] (-0.1,-0.3) -- (24.2,6) node [black, right] {$\frac{q}{n+1}$} ; 
 \draw[gray] (-0.2,-3) -- (24.2,1) node [black, right] {$\frac{q}{w+1}$} ;
 
\draw[black, dashed] (0,2) -- (0,-4) node [black, below] {$q_k$} ; 
\draw[black, dashed] (24,9) -- (24,-4) node [black, below] {$q_{k+1}$} ; 
\draw[black, dashed] (18,9) -- (18,-4) node [black, below] {$p_k$} ; 
\draw[black, dashed] (4,2) -- (4,-4) node [black, below] {$r_k$} ; 
\draw[black, dashed] (10,5) -- (10,-4) node [black, below] {$t_k$} ; 
 \draw[black, dashed] (16,7) -- (16,-4) node [black, below] {$u_k$};
 \draw[black, dashed] (7,5) -- (7,-4) node [black, below] {$s_k^{M}$} ; 
 \draw[black, dashed] (5,1) -- (5,-4) node [black, below] {$s_k^{m}$} ;    
 \draw[black, dashed] (6,1) -- (6,-4) node [black, below] {$s_k$} ;  

 \end{tikzpicture}
 \end{center}
 \caption{Generic system $\bP_{k}^{\delta}$, dotted are the extremal cases $\delta=0$ and $1$.}\label{Fig}
 \end{figure}

Choose a parameter $\delta \in (0,1)$ that will provide uncountability. We first construct elementary systems $\bP_{k}^{\delta}$ on intervals $[q_k,q_{k+1}]$. Denote $\alpha : = -\frac{\log (\eps)}{n+1} + 2R_n$ and $\beta_k := (-\frac{\log (\nu)}{w+1} + 2R_n)$ or $\log q_k$. Figure \ref{Fig} illustrates the construction.\\

At $q_k$, we fix $P_1(q_k) = \cdots = P_n(q_k) = \frac{q_k}{n+1}-\alpha$ and $P_{n+1}(q_k) = \frac{q_k}{n+1}+n\alpha$.
 
Define 
\begin{align*}
p_k &= q_k\frac{w+1}{n+1} - (w+1)(\alpha - \beta_k), &q_{k+1} &= \frac{w}{n} q_k +  \frac{(w-1)(n+1)}{n}(\alpha - \beta_k), \\
 r_k &= q_k + (n^2-1) \alpha,  &u_k &= p_k  - (n+1) \alpha,,\\
s_k^M &= r_k + n \log q_k, &s_k^m &= r_k + \log q_k, \\
s_k &= \delta s_k^m + (1-\delta) s_k^M \in [s_k^m, s_k^M]&t_k &= s_k + (n-1)(s_k-r_k),
\end{align*}
and note that $q_{k+1}$ does not depend on $\delta$.\\

On the interval $[q_k,r_k]$, the $n-1$ components $P_2= \cdots = P_n$ coincide and have slope $\frac{1}{n-1}$ while $P_1$ and $P_{n+1}$ are constant. By definition of $r_k$, $P_2(r_k) = P_{n+1}(r_k)$.\\
On the interval $[r_k,s_k]$, $P_{n+1}$ has slope $1$ while all other component $P_i$ are constant and on $[s_k,t_k]$ the $n-1$ components $P_2 = \cdots = P_{n}$ coincide and have slope $\frac{1}{n-1}$ while $P_1$ and $P_{n+1}$ are constant. By definition of $t_k$, $P_{n+1}(t_k) = P_2(t_k)$.\\
On the interval $[t_k,u_k]$, the $n$ components $P_2= \cdots = P_{n+1}$ coincide and have slope $1/n$, while $P_1$ is constant.\\
On the interval $[u_k,p_k]$, $P_{n+1}$ has slope $1$ while all other component $P_i$ are constant. By definition of $u_k$ and $p_k$, at $p_k$ we have $P_{n+1}(p_k)-P_n(p_k) = (n+1)\alpha$ and $P_1(p_k) = \frac{p_k}{w+1} - \beta_k$. \\
On the interval $[p_k,q_{k+1}]$, $P_1$ has slope $1$ while all other components are constant. By definition of $q_{k+1}$ we have $P_1(q_{k+1}) = \cdots = P_n(q_{k+1}) = \frac{q_{k+1}}{n+1}-\alpha$ and $P_{n+1}(q_{k+1}) = \frac{q_k}{n+1}+n\alpha$.


Choose $q_1$ large enough so that $t_1 <u_1$. Since $w>n$, the inductive sequence $(q_k)_{k\ge1}$ tends to infinity and we define $\bP_{}^{\delta}= \bigcup_{k\ge1} \bP_{k}^{\delta}$ the concatenation of the systems $\bP_{k}^{\delta}$. It is a system, as properties (i) - (iii) holds, in particular at $q_k$.\\
By construction, 
\begin{eqnarray*}
\liminf_{q\to \infty} \frac{P_1^{\delta}(q)}{q} &=&  \liminf_{k\to \infty} \frac{P_1^{\delta}(p_k)}{p_k}  = \frac{1}{1+ w}\\
\min_{[q_k,q_{k+1}]}\frac{q}{n+1} - P_1^{\delta}(q) &=& \alpha  =  -\frac{\log (\eps)}{n+1} + 2R_n\\
\limsup_{q\to\infty} \left(\frac{q}{w+1} - P_{1}^{\delta}(q)\right) &=& \limsup_{k\to\infty} \left(\frac{p_k}{w+1} - P_{1}^{\delta}(p_k)\right)  = \limsup_{k\to\infty} \beta_k 
 \end{eqnarray*}
and the requested properties hold. Furthermore, $P_{1}^\delta(s_k^M) = P_{n+1}(q_k) + \delta \log q_k$, so that
$$\| \bP^{\delta} - \bP^{\delta'} \| \ge |\bP^{\delta}(s_k^M)  - \bP^{\delta'}(s_k^M) | = |\delta - \delta'|\log q_k.$$
Hence, for $\delta \neq \delta'$ the systems $\bP^{\delta}$ and $ \bP^{\delta'} $ are non-equivalent. Uncountability follows. \qed\\

For simultaneous approximation, we rather use the setting and notation of Schmidt and Summerer for parametric geometry of numbers, studying the convex body 
\[ \cC_{\bx}'(q) := \left\{ \by \in \bR^{m+1} \mid |y_1| \le e^{mq} , \max_{1< i \le n} |y_1x_i - y_{i+1}| \le e^{-q}  \right\}\]

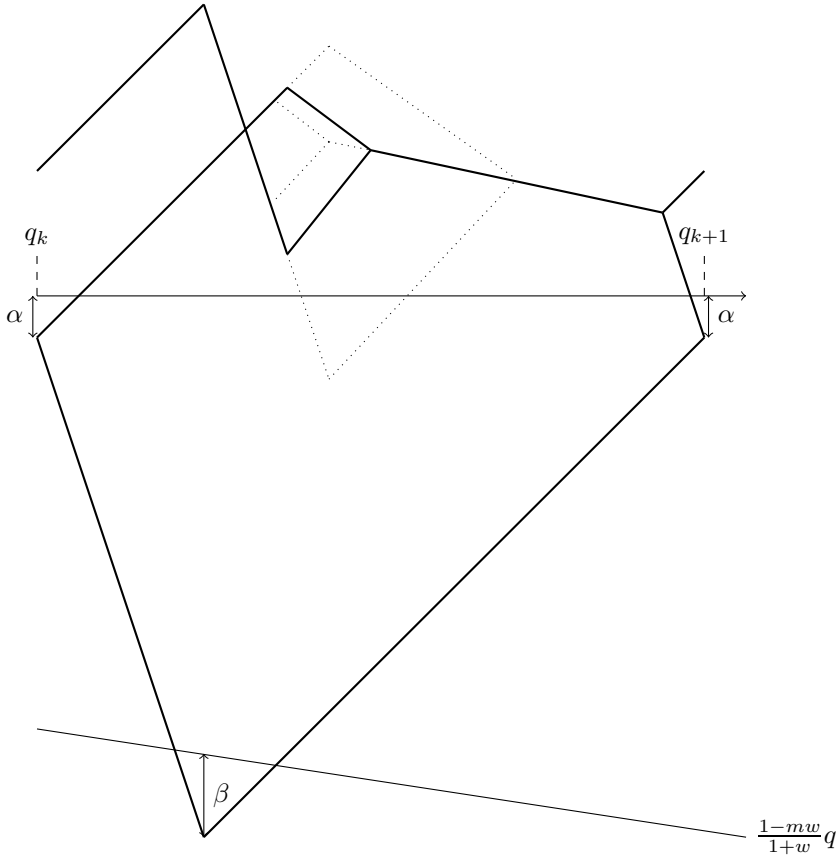
\begin{figure}[h!]
 \begin{center}
 \begin{tikzpicture}[scale=0.55]
 
   \draw[black, ->] (0,0) -- (17,0); 
  \draw[black] (0,-10.4) -- (17, -13) node [black, right] {$\frac{1-mw}{1+w}q$} ;
   \draw[black, <->] (-0.1,0) -- (-0.1,-1) node [black, midway, left] {$\alpha$};
 \draw[black, <->] (16.1,0) -- (16.1,-1) node [black, midway, right] {$\alpha$};
  \draw[black, <->] (4,-13) -- (4,-11) node [black, midway, right] {$\beta$};
 
\draw[black, dashed] (0,0) -- (0,1) node [black, above] {$q_k$} ; 
\draw[black, dashed] (16,0) -- (16,1) node [black, above] {$q_{k+1}$} ; 
 
  \draw[black, thick] (0,3) -- (4,7) ;
  \draw[black, thick] (0,-1) -- (6,5) ;
  \draw[black, thick] (6,1) -- (4,7) ;
  \draw[black, thick] (6,1) -- (8,3.5) ;
  \draw[black, thick] (8,3.5) -- (6,5) ;
  
   \draw[black, thick] (8,3.5) -- (15,2) ;
  
  \draw[black, thick] (0,-1) -- (4,-13) ;
  \draw[black, thick] (16,-1) -- (4,-13) ;
  \draw[black, thick] (16,-1) -- (15,2) ; 
   \draw[black, thick] (15,2) -- (16,3) ;
   
   \draw[black,dotted] (7,3.7) -- (8,3.5);
   \draw[black,dotted] (7,3.7) -- (5.7,2.3);
   \draw[black,dotted] (7,3.7) -- (5.7,4.7);
   
   \draw[black,dotted] (6,1) -- (7,-2);
    \draw[black,dotted] (6,5) -- (7,6);
    \draw[black,dotted] (7,-2) -- (11.5,2.8);
   \draw[black,dotted] (7,6) -- (11.5,2.8);

 \end{tikzpicture}
 \end{center}
 \caption{Generic system for simultaneous approximation.}\label{FigSim}
 \end{figure}

See for example \cite{SS1,SS2,S2020,S2022}. The analogous proof relies on the construction of dual systems depicted by Figure \ref{FigSim}.

\end{document}